\documentclass[review]{elsarticle}

\usepackage[utf8]{inputenc} 
\usepackage{amssymb}
\usepackage{amsmath}
\usepackage{amsfonts}
\usepackage{amsthm}
\usepackage{multirow}
\usepackage{graphicx}
\usepackage{lineno}
\modulolinenumbers[5]

\newcommand{\dparttheta}[1]{ D_\theta^1 \left( #1 \right) }
\newcommand{\dpartthetadois}[1]{ D_\theta^2 \left( #1 \right)  }
\newcommand{\argmin}[1]{\underset{#1}{\operatorname{argmin}}\,}

\newcommand{\logpac}[1]{\log \left(  #1 \right)}
\newcommand{\loglogn}{\log\log n}
\newcommand{\cas}[0]{\underset{a.s.}{\rightarrow}}

\newtheorem{thm}{Theorem}
\newtheorem{cor}{Corollary}
\newtheorem{lem}{Lemma}
\newtheorem{defn}{Definition}

\newcommand{\cupk}{\bigcup\limits_{k=0}^\infty}
\newcommand{\cupkN}{\bigcup\limits_{k\in{\N^q}}}
\newcommand{\numbersys}[1]{\ensuremath{\mathbb{#1}}}

\newcommand{\reals}{\numbersys{R}}

\newcommand{\N}{\numbersys{N}}
\newcommand{\X}{\numbersys{X}}
\newcommand{\M}{\numbersys{M}}
\newcommand{\Y}{\numbersys{Y}}
\newcommand{\I}{\numbersys{I}}
\newcommand{\Prob}{\mbox{\normalfont P}}
\renewcommand{\P}{\Prob}

\newcommand{\go}{\rightarrow}

\newcommand{\thetahat}{ \hat{\theta}}

\newcommand{\alphahat}{ \hat{\alpha}}

\newcommand{\fcalt}{ {\cal F}_{t-1}}

\newcommand{\limitstUKbarN}{ \limits_{t=1+\bar{k}}^n }
\newcommand{\limitstUKbN}{ \limits_{t=1+\bar{k}}^n }
\newcommand{\D}[2]{ \frac{\partial #1}{\partial \alpha_{#2}} }

\newcommand{\DD}[3]{ \frac{\partial^2 #1}{\partial \alpha_{#2} \alpha_{#3}} }
\newcommand{\DDD}[4]{ \frac{\partial^3 #1}{\partial \alpha_{#2} \alpha_{#3} \alpha_{#4}} }

\newcommand{\kbar}{ { \bar{k} } }
\newcommand{\limn}{ \lim\limits_{n \go \infty }}
\newcommand{\limt}{ \lim\limits_{t \go \infty }}
\newcommand{\limsupn}{ \limsup\limits_{n \go \infty }}
\newcommand{\liminfn}{ \liminf\limits_{n \go \infty }}
\newcommand{\logl}[1]{\log L_{n,k}( #1 )}
\newcommand{\loglnk}[1]{\log L_{n,k}( #1 )}
\newcommand{\loglnr}[1]{\log L_{n,r}( #1 )}
\newcommand{\loglnp}[1]{\log L_{n,p}( #1 )}
\newcommand{\loglnl}[1]{\log L_{n,l}( #1 )}

\newcommand{\logLnk}{\log L_{n,k}}

\newcommand{\R}{\numbersys{R}}

\newcommand{\EDC}{\mbox{\normalfont EDC}}
\newcommand{\Tr}{\mbox{\normalfont Tr}}

\newcommand{\norm}[1]{\left\lVert#1\right\rVert}
\newcommand{\Ht}{H_t}

\newcommand{\Hdot}{\dot{H}_t}

\begin{document}
\begin{frontmatter}
\title{Model identification using the Efficient Determination Criterion}
\author[address]{Paulo Angelo Alves RESENDE\corref{mycorrespondingauthor}}
\cortext[mycorrespondingauthor]{Corresponding author}
\ead{pa@pauloangelo.com}
\author[address]{Chang Chung Yu DOREA}
\ead{chang@mat.unb.br}
\address[address]{Department of Mathematics, University of Brasília, Brasília 70910-900, Distrito Federal, Brazil}

\begin{abstract}
In the realm of the model selection context,
Akaike's and Schwarz's information criteria, AIC and BIC, 
have been applied successfully for decades for
model order identification.
The Efficient Determination Criterion (EDC)
is a generalization of these criteria,
proposed originally to define a strongly consistent class
of estimators for the dependency order of a
multiple Markov chain.
In this work, the EDC is generalized 
to partially nested models, which encompass 
many other order identification problems.
Based on some assumptions, 
a class of strongly consistent estimators
is established in this general environment.
This framework is applied to BEKK multivariate GARCH models
and, in particular, the strong consistency of 
the order estimator based on BIC is established for these models.
\end{abstract}

\end{frontmatter}

\linenumbers

\section{Introduction}

The order identification problem was initially dealt by
using nested hypothesis tests
in evaluating the order of multiple Markov chains
\cite{bartlett_frequency_1951,hoel_test_1954,good_likelihood_1955,billingsley_statistical_1961,anderson_statistical_1957},
Autoregressive models
\cite{quenouille_1947,whittle_1951,whittle_1954,bartlett_1953}, among others.
In the selection model context,
Akaike \cite{akaike_new_1974} proposed the use of 
the information criterion AIC, 
aiming to avoid empirical analysis on the estimation process.
Subsequently, Schwarz \cite{schwarz_estimating_1978}
proposed the information criterion BIC.
Since then, these criteria have been applied in such contexts as
selecting models in Autoregressive (AR) and Autoregressive Moving Average (ARMA) process 
\cite{akaike_1979,hannan,hannan_1980,ozaki_1977,ogata_1980,choi_1992},
estimating dependency order in multiple Markov chains
\cite{tong_determination_1975,katz_criteria_1981,csiszar_consistency_2000}, 
detecting change-points in non-homogeneous Markov chains \cite{polansky_2007},
estimating the length of the hidden state space of a hidden Markov model \cite{finesso_1990},
estimating order in 
Autoregressive Conditional Heteroskedasticity process (ARCH) \cite{hughes_selecting_2004}
and on estimating dependency order in specific situations \cite{raftery_model_1985}.

Zhao et al. \cite{zhao_determination_2001},
on estimating the order of multiple Markov chains,
introduced the Efficient Determination Criterion (EDC),
which allows for adjustments on the penalty term used
in the criteria AIC and BIC.
Also, a class of strongly consistent estimators was established in the same work.
Afterwards, Dorea \cite{dorea_optimal_2008} extended this class and proposed the asymptotic optimal order estimator,
which had its better performance verified by the extensive use of numerical simulations \cite{resende_2009}.

In this work, the concept of ``nested models'' is generalized to
class of partially nested models and the EDC criterion
is extended to this new context.
Some results regarding the consistency of EDC order estimators
are established based essentially on assumptions about the likelihood function.
This approach is applied to 
state the consistency of the BIC order estimator for
BEKK multivariate GARCH models, 
which encompass the univariate version GARCH as particular case.

Section 2 provides the general results, that may be applied
in a variety of models to establish the EDC order estimators.
Section 3 presents the approach applied to 
BEKK multivariate GARCH models. 
The proofs of the stated results are in 
the  appendices.

\section{General framework}

The essence of nested models have being used since 
the pioneer researches using hypothesis tests. 
However, practically all works focused on particular
cases and the formal definition and treatment  
of the concept of nested models were unused.
Nishii \cite{nishii_maximum_1988} 
firstly proposed a general estimator 
for the dimension of i.i.d. models.
A relevant piece of Nishii's technique 
is adapted to our purposes.

For an arbitrary time discrete stochastic process $\X=\left\{ X_t \right\}_{t\in \N }$,
$E \subseteq \reals^p$ the set of possible values of $X_t$ and
$\nu$ a fixed measure on $E$,
we define a family of statistical models for $\X$ as
     $$
        M=\left\{ f(x_1^n,\theta,n) : \; \theta \in \Theta,\; n\geq 1 \right\}
     $$
{\raggedright where $f(x_1^n,\theta,n)$ represents the set of possible densities for $x_1^n$
with respect to the product measure on $E^n$,
which depends on the parameter $\theta \in \Theta \subseteq \reals^d$, 
and $x_1^n = x_1 x_2 \dots x_n$ is a realization of $\X$.}
We may denote $f(x_1^n,\theta)=f(x_1^n,\theta,n)$ to simplify the notation.

Two statistical models
     $$
        M_k=\left\{ f(x_1^n,\theta,n) : \; \theta \in \Theta_k,\; n\geq 1 \right\}
        \;\;\;\;and
     $$
     $$
        M_p=\left\{ f(x_1^n,\theta,n) : \; \theta \in \Theta_p,\; n\geq 1 \right\}
     $$
are nested, denoted by $M_k\subseteq M_p$, if $ \Theta_k \subseteq \Theta_p$ 
and, for all $\theta \in \Theta_k$, $x_1^\infty \in E^\infty$, exists $c\in(0,\infty)$
such as
    $$
       \limn \frac{f_k(x_1^n,\theta)}{f_p(x_1^n,\theta)} = c.
    $$
For $q\in\N$, $p=(p_1,\dots,p_q)\in \N^q$ and  $k=(k_1,\dots,k_q)\in\N^q$, 
we define the usual order relation $p\geq k$ iff $p_i \geq k_i$ for $i=1\dots q$, which
makes $(\N,\geq)$ a partially ordered set. 
For $p\not\geq k$ we mean that $p<k$ or $p$ and $k$ are not related.
The set $\M=\left\{ M_k \right\}_{k\in \N^q}$ is a class of partially nested models
if
$$
  M_k \subseteq M_p  \Leftrightarrow k \leq p.
$$
We say that an element $m_r\in \cupkN M_k$ has order $r\in\N^q$ if $m_r\in M_r$ and
$m_r\in M_k$ implies that $M_r\subseteq M_k$. 
In this context, for a sample $x_1^n$ 
and $\thetahat_k$ the maximum likelihood estimator (MLE) of $\theta$ supposing the order $k$,
$\left\{L_{n,k}(x_1^n,\theta)\right\}_{(n,k)\in\N\times \N^q}$
is a class of functions $L_{n,k}:E^n\times \Theta_k \go \R$ that satisfies
\begin{equation}
\label{eq:Lnk01}
L_{n,k}(x_1^n,\thetahat_k)=\sup\limits_{\theta\in\Theta_k} \left\{L_{n,k}(x_1^n,\theta)\right\}
\end{equation}
and for $\theta\in\Theta_k$ and $p\geq k$,
\begin{equation}
\label{eq:Lnk02}
L_{n,p}(x_1^n,\theta) \geq L_{n,k}(x_1^n,\theta) \;\;\;\;and\;\;\;\; \limn \frac{L_{n,p}(x_1^n,\theta)}{L_{n,k}(x_1^n,\theta)} < \infty.
\end{equation}
To simplify notation, we shall denote $L_{n,k}(\theta)=L_{n,k}(x_1^n,\theta)$.
In most situation, the $L_{n,k}$ functions are merely the likelihood for each $n$ and $k$.
Now, we define the EDC estimator for a class of partially nested models.

\begin{defn}
Let $\M$ be a class of partially nested models,
$m_r\in \cupkN M_k$ of order $r$ and $K\geq r$. 
The EDC estimator is defined by
\begin{equation}
 \label{eq:rEDC}
  \hat{r} = \argmin{k\leq K} \left\{ \EDC (k) \right\} 
\end{equation}
for
$$
  \EDC(k) = -\logl{\thetahat_k} + c_n \gamma(k),
$$
$c_n$ a sequence of positive numbers and  $\gamma(k) = \dim(\Theta_k)$.
\end{defn}

We need the following assumptions to conclude consistency 
for the EDC estimator based on the asymptotic behaviour of
the $c_n$ sequence. In what follows, $r$ is the order of $\X$, 
$\theta_r$ is the true parameter, i.e. the one that gives the density for the process $\X$ 
and $\thetahat_k$ is the MLE of $\theta_r$ supposing the order $k$.

\paragraph{\textbf{Assumption A1}} For all $k\geq r$, $\theta_r$ is an interior point of $\Theta_k$ and
      $$
          \thetahat_k \cas \theta_r.
      $$

\paragraph{\textbf{Assumption A2}}  For all $k,n\in\N$, $\logl{x_1^n,\theta}$ and its derivatives\\*
                           $D^1_\theta (\logl{x_1^n,\theta})$,
                           $D^2_\theta (\logl{x_1^n,\theta})$ and
                           $D^3_\theta (\logl{x_1^n,\theta})$
                           are measurable with respect to $x_1^n$ and 
                           continuous with respect to $\theta$.

\paragraph{\textbf{Assumption A3}}
For $k\geq r$, there exists $c<\infty$ and a symmetric and positive definite matrix $A_2$, such as,
for all $\dot{\theta} = (1-s)\thetahat_k + s\theta_r$, $s\in(0,1)$, $(i,j,l)\in \left\{ 1,\dots,\gamma(k) \right\}^3$,
$$
\limn \frac{\left( D_\theta^3 (\logl{\dot{\theta}}) \right)_{i,j,l} }{n} < c  \;\;\;\; a.s. \;\;\;\; and
$$
$$
\limn \frac{D_\theta^2 (\logl{\dot{\theta}}) }{n} = A_2  \;\;\;\; a.s. 
$$

\paragraph{\textbf{Assumption A4}} If $k\geq r$,
           \begin{equation}
             \label{eq:A4}
               \limsupn \norm{ \frac{ D^1_\theta \logl{\theta_r}  }{ \sqrt{2 n \loglogn}}  }  < \infty  \;\;\;\; a.s.
           \end{equation}

\paragraph{\textbf{Assumption A5}}  If $k\not\geq r$,
           $$
              0<\limn \frac{\loglnr{\thetahat_r} - \loglnp{\thetahat_k}}{n} \;\;\;\; a.s.
           $$

We note that Basawa \& Heyde \cite{basawa_1976} propose
basically the use of assumptions A1-A3 and an analogous of A4
to conclude the asymptotically normality of the parameter estimator $\thetahat_k$.
The approach used here is quite similar to their.
The next result establishes the class of consistent EDC order estimators
based on the assumptions and on the asymptotic behaviour of the sequence $c_n$.
The proof is in the Appendix A.
The Corollary concludes the strong consistency for the BIC order estimator
defined in this general context.

\begin{thm}
\label{thm:consistencyEDC}
Let $\X$ be a discrete time stochastic process taking values in $\R^m$,
$\M$ its respective class of partially nested models,
$m_r\in\cupkN M_k$ of order $r$, $\hat{r}$ as defined in \eqref{eq:rEDC},
and the assumptions A1-A5 are satisfied.
Then $\hat{r}\cas r$ if
$$
  \liminfn \frac{c_n}{\loglogn} = \infty \;\;\;\; and \;\;\;\; \limn \frac{c_n}{n} =0.
$$
\end{thm}

\begin{cor} Supposing the same hypothesis of Theorem \ref{thm:consistencyEDC}, the BIC
order estimator bellow is strongly consistent.
$$
  \hat{r}_{bic} = \argmin{k\leq K} \left\{ -\logl{\thetahat_k} + \frac{\log n }{2} \gamma(k) \right\} 
$$
for a known $K\geq r$.
\end{cor}

The proof consists in the determination of the asymptotic behaviour of the
differences bellow for arbitrary $p\not\geq r$ and $k>r$.
\begin{equation}
\label{eq:difLogVeron}
\frac{\loglnr{\thetahat_r} - \loglnp{\thetahat_p}}{n}
\end{equation}
\begin{equation}
\label{eq:difLogVerologlogn}
\frac{\loglnk{\thetahat_k} - \loglnr{\thetahat_r}}{\loglogn}
\end{equation}

Assumption A5 is precisely \eqref{eq:difLogVeron}.
For \eqref{eq:difLogVerologlogn} we use assumptions A1 and A2
to enable the use of Taylor series and state 
\begin{scriptsize}
\begin{eqnarray*}
  \limsupn  \frac{ \loglnk{\thetahat_k} - \loglnk{\theta_r} }{\log\log n}   
      & \leq & \limsupn   \frac{  \sqrt{n} (\theta_r-\thetahat_k) }
                            { \sqrt{2 \loglogn}}
                       \frac{- \dpartthetadois{\loglnk{\thetahat_k}} }
                            {n}
                       \frac{ \left( \sqrt{n} (\theta_r-\thetahat_k)\right)'}
                            { \sqrt{2 \loglogn}} \notag \\
      &  & + o(\loglogn)
\end{eqnarray*}
\end{scriptsize}
and
\begin{eqnarray*}
   \frac{
     1
   }{\sqrt{n \loglogn}} 
   D^1_\theta \loglnk{\theta_r} 
   \left( \frac{- \dpartthetadois{\loglnk{\thetahat_k}} } {n} \right)^{-1} 
   & = &  \frac{\sqrt{n}}{\sqrt{\loglogn}} (\thetahat_k - \theta_r). 
\end{eqnarray*}
Using A3 and A4 we define a upper bound for \eqref{eq:difLogVerologlogn}.

Defining the asymptotic behaviour of $D_\theta^1 \logLnk(\theta_r)$ is generally
easier when compared to the effort in manipulating directly the equation 
\eqref{eq:difLogVerologlogn} to state its asymptotic behaviour.
Besides that, assumptions A1-A3 are commonly used to establish the asymptotic
normality for the MLE and usually are available in the literature.

\section{BEKK-GARCH order estimation}

Engle \cite{engle_1982} originally proposed 
the use of ARCH models on modelling time series in economy.
His work has been hugely influential in the area and motivated many
generalizations and/or adaptations such as 
GARCH \cite{bollerslev_generalized_1986}, NGARCH \cite{engle_1993}, EGARCH \cite{nelson_1991} 
and the multivariate generalizations BEKK-GARCH \cite{engle_1995}, 
VEC-GARCH \cite{bollerslev_1988}, CCC-GARCH \cite{bollerslev_1990}, 
Factor-GARCH \cite{engle_1990}, among others.

The multivariate models have special applications in portfolio selection and asset pricing.
In this family, the BEKK-GARCH models has particular relevance due to its generality and
the amount of research available in the literature. Among the mentioned, only the
VEC-GARCH is more general than the BEKK-GARCH model. However, the VEC-GARCH cases that can not
be represented in the BEKK-GARCH parametrization are somewhat degenerated \cite{boussama_2011,stelzer_2008}.

Boussama \cite{boussama} 
immersed the BEKK-GARCH models into the framework of general state space Markov chains
and used algebraic topology to conclude the geometric ergodicity of such models under
regularity conditions. This work is also published in \cite{boussama_2011} with minor changes.

Comte \& Lieberman \cite{comte_2003} used Boussama's results to prove the consistency conditions
proposed by Jeantheau \cite{jeantheau_1998} and the conditions proposed by
Basawa \& Heyde \cite{basawa_1976} and conclude the strong consistency and asymptotic
normality for the MLE of the parameter $\theta_r$.

As with the particular case of ARCH models, 
until now, there is no results regarding consistency of order
estimators for BEKK-GARCH models. However, the AIC and BIC
information criteria have been used without further formalization.
In what follows, we present some preliminary results, which are
used to prove the assumptions A1-A5 and conclude 
the consistency of the EDC order estimator $\hat{r}$ for such models,
which encompass the consistency of the BIC order estimator as a particular case.

For $k=(k_1,k_2)\in\N^2$, a random sequence $\X=\left\{ X_t \right\}_{t\in \N}$
taking values in $\R^m$ is a BEKK-GARCH(k) model  if for all $t\in\N$, 
\begin{equation}
\label{eq:relacaoBEKK}
                 X_t = \left( H_t \right)^{\frac{1}{2}} \varepsilon_t ,
\end{equation}
where, for $m\times m$ matrices $C$, $\left\{A_{ls}\right\}$ and $\left\{B_{ls}\right\}$, $C$ positive definite and $N\in\N$,
              \begin{equation*}
                  H_t = C + 
                            \sum\limits_{l=1}^{k_2} 
                            \left( 
                            \sum\limits_{s=1}^N 
                                A_{ls} X_{t-l} X'_{t-l} A'_{ls}
                            \right)
                            +
                            \sum\limits_{l=1}^{k_1}
                            \left( 
                            \sum\limits_{s=1}^N 
                                B_{ls} H_{t-l} B'_{ls}
                            \right),
               \end{equation*}
$\left\{ \varepsilon_t \right\}_{t\in\N} \sim {\cal N}(0,I_m)$, and $I_m$ is the $m\times m$ identity matrix.
The process $\X$ can be represented as a Markov chain $\Y=\left\{ Y_t \right\}_{t\in\N}$ defined by
   \begin{equation*}
 Y_t = (vech(H_{t+1})',vech(H_t)',\dots,vech(H_{t-k_1+2})',X'_t,X'_{t-1},\dots,X'_{t-k_2+1})',
   \end{equation*}
where $vech$ is the operator that stacks the lower triangular portion of a matrix.
Boussama proved that $\Y$ is a positive Harris and geometric ergodic Markov chain if
   \begin{equation}
     \label{eq:ergodicidadeBEKKGARCH}
         \rho\left( \sum\limits_{l=1}^{k_2} \tilde{A}_l +  \sum\limits_{l=1}^{k_1} \tilde{B}_l \right) < 1,
   \end{equation}
where $\rho$ is the spectral radius,
  \begin{equation*}
     \tilde{A}_l =  D_m^+ \sum\limits_{s=1}^N (A_{ls} \otimes A_{ls}) D_m
     \;\;\;\;and\;\;\;\;
     \tilde{B}_l =  D_m^+ \sum\limits_{s=1}^N (B_{ls} \otimes B_{ls}) D_m,
  \end{equation*}
$vec$ is the operator that stacks the columns of a matrix,
$\otimes$ is the Kronecker product, $D_m$ and $D_m^+$ 
are defined by the operators that satisfy
\begin{equation*}
  vec(A)=D_m vech(A) \;\;\;\;and
\end{equation*}
\begin{equation*}
  vech(A)=D_m^+ vec(A).
\end{equation*}

For $k=(k_1,k_2)$ and a fixed $N\geq 1$, the BEKK-GARCH(k) model can be immersed 
in a class of partially ordered 
models considering, for $p>k$ and $\theta_k\in\Theta_k$,
$\Theta_k\subsetneq \Theta_p$, 
the matrices
$\left\{A_{ls}(\theta_k)\right\}$, 
$\left\{B_{l's}(\theta_k)\right\}$ and
$C(\theta_k)$ for $l=1\cdots k_2$, $l'=1\cdots k_1$ and $s=1\cdots N$.
Denoting $\kbar=\max(k_1,k_2)$, we have
   \begin{equation*}
      f(x_1^n,\theta_{k}) = C_1(x_1^\kbar) \prod\limitstUKbarN \frac{1}{\sqrt{(2\pi)^{m/2} \det(H_{t,\theta_k})}} \exp\left( -\frac{1}{2 } x'_t H_{t,\theta_k}^{-1} x_t \right)
   \end{equation*}
for
               \begin{equation*}
                  H_{t,\theta_k} = C + 
                            \sum\limits_{l=1}^{k_2} 
                                A_{l} X_{t-l} X'_{t-s} A'_{l}
                            +
                            \sum\limits_{l=1}^{k_1}
                                B_{l} H_{t-l, \theta_k} B'_{l}.
               \end{equation*}
Also, the following definition for the functions $\logLnk$ satisfies \eqref{eq:Lnk01} and  \eqref{eq:Lnk02}.
               \begin{equation*}
               \logl{\theta_k} = \sum\limits_{1+\bar{k}}^n l_t(\theta_k)
               \end{equation*}
               for
               \begin{equation}
               \label{eq:logveroBEKKGARCH}
               l_t(\theta_k) = -\frac{1}{2} X'_t H^{-1}_{t,\theta_k} X_t -\frac{1}{2} \log \det(H_{t,\theta_k}).
               \end{equation}

The nesting relation $M_k\subset M_p$ can be observed taking 
$A_{l s}$ and $B_{l' s}$ as null matrices for $l>k_2$ and $l'>k_1$.

In particular, if $N=1$ and $\theta_k$ is the columns of the matrices 
$\left\{A_{l}\right\}$, 
$\left\{B_{l'}\right\}$ and
$C$, we may construct
$\Theta_k\subseteq \R^{m^2(2\kbar +1)}$, 
   \begin{equation*}
       \Theta_k = \Omega_0 \times \Omega_1\times \dots \times \Omega_{2\kbar},
   \end{equation*}
   for $\Omega_i =\{0\}^{m^2}$ if $i/2> k_2$ and $i$ is odd or  $i/2> k_1$ and $i$ is even, 
   for the remaining cases, $\Omega_i \subseteq \R^{m^2}$ has non-empty interior.
   Assuming $A_i=0$ if $i > k_2$, $B_i=0$ if $i > k_2$, 
   \begin{equation*}
      \theta_{k}=(vec(C),vec(A_1),vec(B_1),\dots,vec(A_{\bar{k}}),vec(B_{\bar{k}}))\in\Theta_{k}.
   \end{equation*}
   In this case,
   \begin{equation*}
     \gamma(k) =  m^2 ( 1+ k_1 + k_2).
   \end{equation*}

For the order, $r$, of a BEKK-GARCH process $\X$, we consider the lowest $k$ such as
$\X$ can be represented by \eqref{eq:relacaoBEKK}.
Assuming the following conditions (B1-B5), we establish the
Theorem \ref{thm:bekk} that concludes assumptions A1-A5 and states the
class of strong consistent EDC order estimators.
For $k\geq r$,

\paragraph{\textbf{B1}} $\Theta_k$ is compact and $\theta_r$ is an interior point of $\Theta_k$.
\paragraph{\textbf{B2}} There exists a $c>0$ such as $\inf_{\theta\in\Theta_k} \det C(\theta)>c$.
\paragraph{\textbf{B3}} The model is identifiable, i.e. $H_t(\theta) = H_t(\theta')$ a.s. if and only if $\theta=\theta'$.
\paragraph{\textbf{B4}} $C(\theta)$, $\tilde{A}_l(\theta)$ and $\tilde{B}_l(\theta)$ 
                        and their derivatives, with respect to $\theta$, until order 3 are continuous.
\paragraph{\textbf{B5}} $X_t$ admits bounded moments of order 16.

Comte \& Lieberman use B1-B4 and finite moments of order 8 in B5
to conclude the asymptotic normality of $\thetahat_k$.
We need the finiteness for moments of order 16 to conclude assumption A4.

\begin{thm}
\label{thm:bekk}
Let $\X$ be a BEKK-GARCH(r) of order $r$, satisfying 
\eqref{eq:ergodicidadeBEKKGARCH} 
and conditions B1-B5. Then the EDC order estimator defined at \eqref{eq:rEDC}
is strongly consistent if
$$
  \liminfn \frac{c_n}{\loglogn} = \infty \;\;\;\; and \;\;\;\; \limn \frac{c_n}{n} =0.
$$
\end{thm}

Boussama \cite{boussama} concluded the geometric ergodicity for the 
associated Markov chain of $\X$ and, in particular,
enabled the use of the Strong Law of Large Numbers (SLLN), which can
be found in Meyn \& Tweedie \cite{meyn_markov_1993}.
However, the geometric ergodicity is not sufficient
to conclude the Law of Iterated Logarithm (LIL), needed to
prove assumption A4.
To overcome this, we use the LIL bellow, stated for square integrable
Martingales, which can be found in \cite{hall_martingale_1980}.
Also, some auxiliary results stated by Comte \& Lieberman \cite{comte_2003}
are used to conclude assumptions A1-A5.

\begin{thm}[Hall \& Heyde (1980)]
\label{thm:LILMartingale}
Let $\left\{ S_n, \fcalt \right\}$ be a martingale, $S_n=\sum_{t=1}^n U_t$, $E(S_n)=0$, $E(S_n^2)<\infty$, 
$\left\{ Z_t \right\}_{t\in\N}$ and
$\left\{ W_n \right\}_{n\in\N}$ 
non-negative random variables 
such as $Z_t$ and $W_t$ are
$\fcalt$ measurable. If
\begin{equation}
\label{thm:LILMartingale.L1}
   \limn \frac{
                \sum\limits_{t=1}^n  U_t \I(|U_t|>Z_t)  - E[U_t \I(|U_t| > Z_t)|\fcalt]
              }{  
                  \sqrt{2 W_n^2 \log\log W_n^2 } 
              } = 0 \;\;\;\;a.s.,
   \tag{L1}
\end{equation}
\begin{equation}
\label{thm:LILMartingale.L2}
   \limn \frac{
                \sum\limits_{t=1}^n  E[U_t^2 \I(|U_t|\leq Z_t)|\fcalt ] -  E[U_t \I(|U_t|\leq Z_t)|\fcalt ]^2
              }{  
                  W_n^2
              } = 1 \;\;\;\;a.s.,
   \tag{L2}
\end{equation}
\begin{equation}
\label{thm:LILMartingale.L3}
   \limn \sum\limits_{t=1}^n 
          \frac{
                E[U_t^4 \I(|U_t| \leq Z_t)| \fcalt]
              }{  
                  W_t^4 
              } < \infty \;\;\;\;a.s.,
   \tag{L3}
\end{equation}
\begin{equation}
\label{thm:LILMartingale.L4}
   \limn \frac{
                W_n
              }{  
                  W_{n+1}
              } = 1 \;\;\;\;a.s.\;\;and\;\;\;\;
    \limn W_n = \infty\;\;\;\;a.s.
   \tag{L4}
\end{equation}
{\raggedright Then}
\begin{equation*}
   \limsupn \frac{S_n}{ \sqrt{2 W_n^2 \log\log W_n^2 } } = 1 \;\;\;\;a.s.
\end{equation*}
{\raggedright and}
\begin{equation*}
   \liminfn \frac{S_n}{ \sqrt{2 W_n^2 \log\log W_n^2 } } = -1 \;\;\;\;a.s.
\end{equation*}
\end{thm}

\section{Conclusion}

The Efficient Determination Criterion (EDC) raises as a promising
approach in the context of partially nested models. 
Mainly because the assumptions A1-A5 simplify the
establishment of strongly consistent order estimators in a variety of models.
Some of these assumptions, for each case, can be found in the literature 
on defining the asymptotic normality for the respective MLE.

Hafner \& Preminger \cite{hafner_2009b} state some results for VEC-GARCH(1,1) models.
If it is possible to generalize these results for arbitrary $k\in \N^2$,
the EDC order estimator can be easily defined for VEC-GARCH models.

As future works, we suggest to weaken the hypothesis B5 and
state the consistency of the EDC estimator for $c_n = O(\loglogn)$.

\section*{Acknowledgement}
The authors gratefully acknowledge the attention of Dr. Farid Boussama in sending his thesis 
and indirectly encouraging the establishment of the EDC for BEKK-GARCH models.

\appendix
\section{Proof of Theorem 1}

The following Lemma is an adaptation of results that can be found at \cite{dorea_optimal_2008}.

\begin{lem}
\label{lem:casosEDC}
    Let $\X$ be a discrete time stochastic process with values in $\R^m$, $\M$ its respective class of partially nested
    models, $m_r\in \cupk M_k$ of order $r$ and $\hat{r}$ as defined in \eqref{eq:rEDC}. 
    Then $\hat{r}$ is strongly consistent ($\hat{r}\cas r$) if,
      for $k \not\geq r$, exists $c_1 \in (0,\infty)$ such as
         \begin{equation}
           \label{eq:h1}
             \limn \frac{\loglnr{\thetahat_r} - \loglnk{\thetahat_k}}{n} \geq c_1 \; a.s. ,
         \end{equation}
      for $k > r$, exists $c_2 \in (0,\infty)$ such as
         \begin{equation}
           \label{eq:h2}
             \limsupn \frac{\loglnk{\thetahat_k} - \loglnr{\thetahat_r}}{\loglogn} \leq c_2 (\gamma(k)-\gamma(r)) \; a.s.
         \end{equation}
    {\raggedright and $c_n$ satisfies}
 \begin{equation}
 \label{thm:convergenciaEDC.e01}
   \limn \frac{c_n}{n} = 0 \;\;\;\;  and \;\;\;\; \liminfn \frac{c_n}{\loglogn} \geq c_2 .
 \end{equation}

\begin{proof}
      We have that
      \begin{eqnarray}
        \left( -\loglnp{\thetahat_p} + \gamma(p) c_n\right)  -  \left( - \loglnl{\thetahat_l} + \gamma(l) c_n \right)  & & \notag \\
             & & \hspace{-7cm} =  \left( \loglnl{\thetahat_l} - \loglnp{\thetahat_p} \right) - c_n \left( \gamma(l) - \gamma(p) \right) .
       \label{thm:convergenciaEDC.p01}
      \end{eqnarray}

{\raggedright Taking $p=r$ and $l=k$ in \eqref{thm:convergenciaEDC.p01} and using \eqref{eq:h2} we get}
\begin{scriptsize}
     \begin{eqnarray*}
       \limsupn \frac{ \left(  - \logl{\thetahat_r} + \gamma(r) c_n  \right) - \left( -\logl{\thetahat_k} + \gamma(k) c_n \right) }{\log \log n}  \\
           & & \hspace{-4cm} \leq  c_2 (\gamma(k)-\gamma(r)) - \liminfn\left( \frac{c_n}{\log\log n} \right) \left( \gamma(k)-\gamma(r) \right) \;\;\;\; a.s. \notag \\
           & & \hspace{-4cm} \leq  c_2 (\gamma(k)-\gamma(r)) - c_2 (\gamma(k)-\gamma(r)) \;\;\;\; a.s.  \notag \\
           & & \hspace{-4cm}  =    0 .
     \end{eqnarray*}
\end{scriptsize}
{\raggedright In the same manner, but taking $l=r$ and $p=k<r$ in \eqref{thm:convergenciaEDC.p01}, and using \eqref{eq:h1} and \eqref{thm:convergenciaEDC.e01}, we have}
 \begin{eqnarray*}
  \liminfn \frac{\left( -\logl{\thetahat_k} + \gamma(k) c_n\right)  -  \left( - \logl{\thetahat_r} + \gamma(r) c_n \right) }{n} \\
       & & \hspace{-4cm} \geq  c_1 - \limsupn \frac{c_n}{n} \left( \gamma(r)-\gamma(k) \right) \;\;\;\;a.s. \notag \\
       & & \hspace{-4cm} >  0 \;\;\;\;a.s.
 \end{eqnarray*}
   Then, using its definition, we conclude that $\hat{r} \cas r$.
\end{proof}
\end{lem}

\begin{lem}
\label{lem:tetaImplicaLog}
    Let $\X$ be a discrete time stochastic process with values in $\R^m$, $r$ its order, $\M$ its respective class of partially nested
    models, the $\logLnk$ functions as defined above and
    $\thetahat_k \in \Theta_k$ the MLE of the true parameter $\theta_r \in \Theta_k$. If assumptions A1-A5 are true, then
 \begin{equation*}
   \limsupn   \frac{ \loglnk{\thetahat_k} - \loglnr{\thetahat_r}}{\log\log n }   \leq \frac{2 c^2}{ \lambda_{\gamma(k)}}   \;\;\;\; a.s.,
 \end{equation*}
 {\raggedright where $\lambda_{\gamma(k)}$ is the lowest eigenvalue of $A_2$ and $c$ is an upper bound for \eqref{eq:A4}.}

\begin{proof}
  Using A1-A2, for large enough $n$, we can take the Taylor expansion of $\loglnk{\theta_r}$ at $\thetahat_k$, which gives
\begin{eqnarray}
  \loglnk{\theta_r} 
    & = & \loglnk{\thetahat_k} + (\theta_r-\thetahat_k) \dparttheta{\loglnk{\thetahat_k} } \notag \\
    & & \hspace{-1cm} 
        + \frac{1}{2}(\theta_r-\thetahat_k) \dpartthetadois{\loglnk{\thetahat_k}} (\theta_r-\thetahat_k)^T 
        + r_n(\theta_r - \thetahat_k)
    \label{lem:tetaImplicaLog.p01}
\end{eqnarray}
   {\raggedright where, for 
                         $\theta_r=(\alpha_1,\cdots,\alpha_{\gamma(k)})$, $\thetahat_k=(\alphahat_1,\cdots,\alphahat_{\gamma(k)})$ and 
                         $\dot{\theta} = (1-s) \theta_k +s \theta_r$, $s\in (0,1)$, }
\begin{eqnarray*}
    r_n(\theta_r - \thetahat_k) 
    & = & \frac{1}{3!} \sum\limits_{i,j,l} \left( D_\theta^3 \loglnk{\dot{\theta}} \right)_{i,j,l} 
    (\alpha_i-\alphahat_i)
    (\alpha_j-\alphahat_j)
    (\alpha_l-\alphahat_l).
\end{eqnarray*}
   {\raggedright 
   By definition, $\thetahat_k$ maximizes $L_{n,k}$, which gives $\dparttheta{\loglnk{\thetahat_k} } =0$. Organizing \eqref{lem:tetaImplicaLog.p01} and dividing by $\loglogn$, we have}
   \begin{scriptsize}
\begin{eqnarray}
  \limsupn  \frac{ \loglnk{\thetahat_k} - \loglnk{\theta_r} }{\log\log n}   
      & \leq & \limsupn   \frac{  \sqrt{n} (\theta_r-\thetahat_k) }
                            { \sqrt{2 \loglogn}}
                       \frac{- \dpartthetadois{\loglnk{\thetahat_k}} }
                            {n}
                       \frac{ \left( \sqrt{n} (\theta_r-\thetahat_k)\right)^T}
                            { \sqrt{2 \loglogn}} \notag \\
      &  & + \limsupn  \frac{|r_n(\theta_r - \thetahat_k) |}{\loglogn}.
    \label{lem:tetaImplicaLog.p03}
\end{eqnarray}
   \end{scriptsize}
   Now, taking the Taylor expansion of $ D_\theta^1 \loglnk{\thetahat_k} $ at $\theta_r$,
   \begin{eqnarray*}
   (0,\cdots,0) & = & D^1_\theta \loglnk{\thetahat_k}  \\
       & = & D^1_\theta \loglnk{\theta_r} + (\thetahat_k - \theta_r)  D^2_\theta \loglnk{\dot{\theta}} ,
   \end{eqnarray*}
   {\raggedright where $\dot{\theta}=s\theta_r +(1-s)\thetahat_k$ and $s\in(0,1)$. Organizing, we have}
   \begin{eqnarray*}
        \frac{1}{\sqrt{n \loglogn}} D^1_\theta \loglnk{\theta_r} & = & \frac{\sqrt{n}}{n \sqrt{\loglogn}} 
                 \left\{ (\thetahat_k - \theta_r) D^2_\theta \loglnk{\dot{\theta}} \right\} \notag \\
       & = &  - \frac{\sqrt{n}}{\sqrt{\loglogn}} (\thetahat_k - \theta_r) \left[ \frac{D^2_\theta \loglnk{\dot{\theta}}}{n} \right] .
   \end{eqnarray*}
   Using that $A_2$ is positive definite we conclude that it is invertible and, for large enough $n$,
   \begin{equation*}
      A_n:= - \left[ \frac{D^2_\theta \loglnk{\dot{\theta}}}{n} \right]
   \end{equation*}
   {\raggedright has inverse $A^{-1}_n$, then }
   \begin{eqnarray}
      \label{lem:tetaImplicaLog.p07a}
        \frac{1}{\sqrt{n \loglogn}} D^1_\theta \loglnk{\theta_r} A_n^{-1} & = &  \frac{\sqrt{n}}{\sqrt{\loglogn}} (\thetahat_k - \theta_r). 
   \end{eqnarray}
{\raggedright Using A4 and \eqref{lem:tetaImplicaLog.p07a}, 
               considering $P_i:\R^{\gamma(k)}\go \R$ as the projection of coordinate $i$, we have}
   \begin{eqnarray*}
       \label{lem:tetaImplicaLog.p05b}
       \limsupn \left| \frac{\sqrt{n} (\alphahat_i - \alpha_i)}{ \sqrt{ 2 \loglogn}} \right| 
       & = & \limsupn  \left| P_i \left(  \frac{\sqrt{n}}{\sqrt{2 \loglogn}} (\thetahat_k - \theta_r)       \right) \right| \notag \\
       & = & \limsupn \left | P_i \left(  \frac{1}{\sqrt{2 n \loglogn}} D^1_\theta \loglnk{\theta_r} A_n^{-1} \right) \right| \notag \\
       & = & \limsupn \left | P_i \left(  \frac{1}{\sqrt{2 n \loglogn}} D^1_\theta \loglnk{\theta_r} A_2^{-1} \right) \right| \notag \\
       & < & \infty,
   \end{eqnarray*}
   {\raggedright which gives, using A1 and A2, }
\begin{scriptsize}
\begin{eqnarray}
    \limsupn \frac{|r_n(\theta_r - \thetahat_k) |}{\loglogn}  & & \notag\\
        & & \hspace{-1cm} \leq \limsupn \frac{1}{3!} \sum\limits_{i,j,l} \left| \frac{\left( D_\theta^3 \loglnk{\dot{\theta}} \right)_{i,j,l}}{n}  \right|
                  \left| \frac{\sqrt{n} ( \alpha_i-\alphahat_i) }{\sqrt{\loglogn}} \right| 
                  \left| \frac{\sqrt{n} ( \alpha_j-\alphahat_j) }{\sqrt{\loglogn}} \right| 
                  |\alpha_l-\alphahat_l| \notag \\
        & & \hspace{-1cm} \leq   c \sum\limits_{l} \limsupn |\alpha_l-\alphahat_l| = 0 \;\;\;\;a.s.
    \label{lem:tetaImplicaLog.p04}
\end{eqnarray}
\end{scriptsize}
Using A3, \eqref{lem:tetaImplicaLog.p03}, \eqref{lem:tetaImplicaLog.p07a}, \eqref{lem:tetaImplicaLog.p04} and that $A_n \cas A_2$, we have
\begin{scriptsize}
\begin{eqnarray}
   \label{lem:tetaImplicaLog.p07}
  \limsupn  \frac{ \loglnk{\thetahat_k} - \loglnk{\theta_r} }{\log\log n}   & & \notag \\
      & & \hspace{-3cm}\leq  \limsupn   
                        \frac{  \left( D^1_\theta \loglnk{\theta_r} \right)A_n^{-1} }
                            { \sqrt{2 n \loglogn}}
                       \frac{- \dpartthetadois{\loglnk{\thetahat_k}} }
                            {n}
                           \left( \frac{  \left( D^1_\theta \loglnk{\theta_r}  \right)A_n^{-1} }
                            { \sqrt{2 n \loglogn}} \right)^T \notag \\
      & & \hspace{-3cm}=  \limsupn   \frac{   D^1_\theta \loglnk{\theta_r} }
                            { \sqrt{2 n \loglogn}}
                            A_2^{-1} A_2 
                           \left( \frac{  D^1_\theta \loglnk{\theta_r} }
                            { \sqrt{2 n \loglogn}} A_2^{-1} \right)^T \notag \\
      & & \hspace{-3cm}=  \limsupn   \frac{   D^1_\theta \loglnk{\theta_r} }
                            { \sqrt{2 n \loglogn}}
                            A_2^{-1}
                           \left( \frac{  D^1_\theta \loglnk{\theta_r} }
                            { \sqrt{2 n \loglogn}} \right)^T \notag \\
      & & \hspace{-3cm}\leq  \frac{1}{\lambda_{\gamma(k)}} \limsupn   \norm{\frac{D_{\theta}^1 \loglnk{\theta_r}}{\sqrt{2 n \loglogn}}}^2 \notag \\
      & & \hspace{-3cm}\leq  \frac{c_5^2}{\lambda_{\gamma(k)}} \;\;\;\;a.s.
\end{eqnarray}
\end{scriptsize}
We used that $A_2$ is symmetric positive definite so is its inverse.
Applying \eqref{lem:tetaImplicaLog.p07} twice, we conclude the proof.
\end{proof}
\end{lem}

Using Assumption A5 and Lemma \ref{lem:tetaImplicaLog} we have \eqref{eq:h1} and \eqref{eq:h2}. The Theorem follows from Lemma \ref{lem:casosEDC}.

\section{Proof of Theorem 2}

\begin{lem}[Comte \& Lieberman (2003)]
\label{lem:resultadosComteBEKKGARCH}
 Let $\X=\left\{ X_t \right\}_{t\in \N}$ be a BEKK-GARCH(k),
 $\theta_r=(\alpha_1,\dots,\alpha_{\gamma(k)})$ its true parameter,
 $\thetahat_k=(\alphahat_1,\dots,\alphahat_{\gamma(k)}) $ the MLE of $\theta_r$. 
 If conditions \eqref{eq:ergodicidadeBEKKGARCH} and B1-B5 are true, then
\begin{itemize}
  \item[(i)]
   \begin{equation*}
      \limn - \frac{D^2_\theta \loglnk{\theta_r}}{n} =  A_2 \;\;\;\; a.s.,
   \end{equation*}
   {\raggedright where}
   \begin{equation}
    \label{eq:defnA2BEKKGARCH}
    A_2 = -E\left(   \frac{\partial^2 l_t(\theta_r)}{\partial\theta \partial\theta' }   \right).
   \end{equation}
  \item[(ii)] $A_2$ is positive definite.
  \item[(iii)]  For all $i,j,l\in \{1,\dots,\gamma(k)\}$,
   \begin{equation*}
     E\left( \sup\limits_{\norm{\theta-\theta_r} \leq \delta} \left|  \DDD{l_t(\theta)}{i}{j}{l}\right| \right) < c(\delta) .
   \end{equation*}
  \item[(iv)]  For all $i\in \{1,\dots,\gamma(k)\}$, $\D{\loglnk{\theta_r}}{i}$ is a square-integrable Martingale.
  \item[(v)]  The MLE $\thetahat_k$ is strongly consistent.
  \item[(vi)] Exists $c\in(0,\infty)$, which does not depend on $t$ or $\theta$, such as
   \begin{equation*}
     \norm{H_t^{-1}} \leq c.
   \end{equation*}
  \item[(vii)] 
   \begin{equation*}
     E\left( \left|  \logpac{det(H_t(\theta_r))}\right| \right) < \infty.
   \end{equation*}
\end{itemize}
\end{lem}

The following Lemma adapts some results from Comte \& Lieberman to our purposes.

\begin{lem}
\label{lem:comportamento_assintotico_derivadasBEKKGARCH}
 Let $\X=\left\{ X_t \right\}_{t\in \N}$ be a BEKK-GARCH(k),
 $\theta_r=(\alpha_1,\dots,\alpha_{\gamma(k)})$ its true parameter,
 $\logLnk$ as defined early, 
 $\thetahat_k \in \Theta_k $ the MLE of $\theta_r$,
 $\dot{\theta}=s\theta_r + (1-s)\thetahat_k$ and $s\in[0,1]$
 and $B_\delta(\theta_r)\subset \Theta_k$ a neighborhood of $\theta_r$.
 If conditions \eqref{eq:ergodicidadeBEKKGARCH} and B1-B5 are true, then

\begin{enumerate}
  \item[(i)] Exists $c\in (0,\infty)$, such as, for all $i,j,l\in\left\{ 1,\dots,\gamma(k) \right\}$,
     \begin{equation*}
        \limsupn 
           \norm{
                \frac{1}{n} \sum\limits_{t=1}^n
                    \sup\limits_{\theta\in B_\delta(\theta_r)} \DDD{l_t(\theta)}{i}{j}{l}
           }      
         \leq c .
     \end{equation*}
  \item[(ii)] 
   \begin{equation*}
      \limn - \frac{D^2_\theta \logl{\dot{\theta}}}{n} =  A_2 \;\;\;\; a.s.
   \end{equation*}
   {\raggedright for $A_2$ as defined in \eqref{eq:defnA2BEKKGARCH}.}
  \item[(iii)] 
   \begin{equation*}
      E\left( \left| \log\left[ det(H_t(\theta_r)) \right] + X'_t H_t^{-1} X_t\right| \right) < \infty.
   \end{equation*}
\end{enumerate}

\begin{proof}
\begin{itemize}
 \item[(i)] Using item (iii) of Lemma \ref{lem:resultadosComteBEKKGARCH} and the 
             Boussama's results, we just apply the SLLN that can be found at \cite{meyn_markov_1993}.
 \item[(ii)] Analogous to the technique used in Lemma 5 of \cite{hafner_2009}, 
            using that $D^2_\theta l_t(\theta)$ and $D^3_{\theta} l_t(\theta)$ are continuous with respect to 
            $\theta$, that $\thetahat_k$ is strongly consistent (Lemma \ref{lem:resultadosComteBEKKGARCH})
            and the mean value Theorem, we have
            \begin{equation*}
             \norm{
                \frac{1}{n} \sum\limits_{t=1}^n \DD{l_t(\dot{\theta})}{i}{j} 
                - \frac{1}{n} \sum\limits_{t=1}^n \DD{l_t(\theta_r)}{i}{j} 
             } \leq
             \sup\limits_{\theta\in B_\delta(\theta_r)} 
             \left\{ \norm{
                 \frac{1}{n} \sum\limits_{t=1}^n 
                 \frac{\partial}{\partial \theta'}
                 \left( 
                 \DD{l_t(\theta)}{i}{j}
                 \right) 
                 }
                 \cdot 
                 \norm{
                   \dot{\theta} - \theta_r
                 }
              \right\} .
            \end{equation*}
            Using item (i) and the strong consistency of $\thetahat$ we conclude the result.
 \item[(iii)] 
     \begin{eqnarray*}
      E\left( \left| \log\left[ det(H_t(\theta_r)) \right] + X'_t H_t^{-1} X_t\right| \right)
      &\leq& 
      E\left( \left| \log\left[ det(H_t(\theta_r)) \right] \right| \right) 
              + E\left( \left|  X'_t H_t^{-1} X_t\right| \right) \\
      &\leq& 
      E\left( \left| \log\left[ det(H_t(\theta_r)) \right] \right| \right) 
              + E\left( \norm{X_t}^2 \right)E\left( \norm{ H^{-1}_t } \right)
     \end{eqnarray*}
     which is bounded by items (vi) and (vii) of Lemma \ref{lem:resultadosComteBEKKGARCH} and by B3.
\end{itemize}
\end{proof}
\end{lem}

\begin{lem} 
\label{lem:derivada1a_assintoticaBEKKGARCH}
 Let $\X=\left\{ X_t \right\}_{t\in \N}$ be a BEKK-GARCH(k),
 $\theta_r=(\alpha_1,\dots,\alpha_{\gamma(k)})$ its true parameter,
 If conditions \eqref{eq:ergodicidadeBEKKGARCH} and B1-B5 are true, thus,
  for all $i\in\left\{ 1,\dots,\gamma(k) \right\}$,
   \begin{equation*}
      \limsup\limits_{n\go \infty} \frac{\D{\loglnk{\theta_r} }{i}}{\sqrt{2 n  \loglogn}} =   E\left( \D{l_1(\theta_r)}{i}^2 \right)^{1/2}  \; \; \; a.s.,
   \end{equation*}
   \begin{equation*}
    \liminf\limits_{n\go \infty} \frac{\D{\loglnk{\theta_r} }{i}}{\sqrt{ 2 n \loglogn}} = -  E\left( \D{l_1(\theta_r)}{i}^2 \right)^{1/2}    \; \; \; a.s.\;\;\;\;and
   \end{equation*}
\begin{equation*}
  \limsupn \frac{\norm{D^1_\theta \loglnk{\theta_r}}}{\sqrt{2\loglogn}} \leq c \;\;\;\;a.s.
\end{equation*}
{\raggedright for $c\in(0,\infty)$.}
   \begin{proof}
   Consider item (iv) of Lemma \ref{lem:resultadosComteBEKKGARCH} and assume $\fcalt=\sigma(X_1,\dots,X_t)$, $Z_t = t^\delta$, $\delta>1$, 
   \begin{equation*}
    U_t =  \D{l_t(\theta_r)}{i} \;\;\;\;and\;\;\;\; W_n = \left[ n E\left( \D{l_t(\theta_r)}{i}^2 \right)  \right]^{1/2},
   \end{equation*}
   {\raggedright where, by \eqref{eq:logveroBEKKGARCH}, }
   \begin{equation*}
                \D{l_t(\theta_r)}{i} =  \frac{1}{2 } \Tr \left( X_t X'_t   H_t^{-1 } \D{H_t}{i} H_t^{-1} 
                     -  H_t^{-1} \D{H_t}{i} \right).
   \end{equation*}
   To apply Theorem \ref{thm:LILMartingale}, we need to prove conditions L1-L4 below.
\begin{itemize}
\item[(L1)]
 By the Chebyshev's inequality, we have
 \begin{eqnarray}
  \label{lem:derivada1a_assintoticaGARCHBEKK.e03}
    P(|U_t| > Z_t) 
    &=&  P\left( 
                 \left| \D{l_t(\theta_r)}{i}\right| > t^\delta
    \right) \notag \\
    &\leq&   \frac{1}{t^{2\delta}} E\left( \D{l_t(\theta_r)}{i}^2  \right).
 \end{eqnarray}
  Using item (iii) of Lemma \ref{lem:resultadosComteBEKKGARCH},
  \begin{equation*}
    \sum\limits_{t=1}^\infty P(|U_t| > Z_t) \leq 
    E\left( \D{l_1(\theta_r)}{i}^2  \right) \sum\limits_{t=1}^\infty  \frac{1}{t^{2\delta}} < \infty.
  \end{equation*}
  By the Borel-Cantelli Lemma, 
  \begin{equation*}
      P\left(\left\{\omega : \I(|U_t|> t^\delta) = 1 \;\;i.o.\right\}\right)
      =
   P\left(\left\{\omega : |U_t|> t^\delta \;\;i.o.\right\}\right) =0 
  \end{equation*}
{\raggedright and thus }
\begin{equation*}
   \limn \frac{
                \sum\limits_{t=1}^n  U_t \I(|U_t|>Z_t)  - E[U_t \I(|U_t| > Z_t)|\fcalt]
              }{  
                  \sqrt{2 W_n^2 \log\log W_n^2 } 
              } = 0 \;\;\;\;a.s.
\end{equation*}

\item[(L2)]
\begin{eqnarray*}
     E(U_t|\fcalt) 
     &=& \frac{1}{2} \Tr\left[ E(X_t X'_t|\fcalt)  H_t^{-1 } \D{H_t}{i} H_t^{-1} -  H_t^{-1} \D{H_t}{i}  \right] \\
     &=& \frac{1}{2} \Tr\left[ H_t H_t^{-1 } \D{H_t}{i} H_t^{-1} -  H_t^{-1} \D{H_t}{i}  \right] \\
     &=& 0 .
\end{eqnarray*}
{\raggedright Using item (iv) of Lemma \ref{lem:resultadosComteBEKKGARCH},}
\begin{eqnarray*}
     E( E(U^2_t|\fcalt) )
     &=& E( U^2_t) < \infty
\end{eqnarray*}
{\raggedright and then, by the SLLN, we have}
\begin{equation*}
   \limn \frac{
                \sum\limits_{t=1}^n  E[U_t^2|\fcalt] 
              }{  
                  n E\left( \D{l_t(\theta_r)}{i}^2  \right) 
              } = 1 \;\;\;\;a.s.
\end{equation*}
{\raggedright By the Dominated Convergence Theorem,}
\begin{equation*}
    \limt E(U_t \I(|U_t| \leq n ) | \fcalt) = 0\;\;\;\;a.s.
\end{equation*}
{\raggedright Considering an arbitrary $\varepsilon>0$, it is required to find a $t$-summable 
               upper bound for}
\begin{equation*}
P\left[|E(U^2_t|\fcalt)-E(U^2_t I(|U_t|\leq t^\delta)|\fcalt)|>\varepsilon\right]
\end{equation*}
{\raggedright and apply the Borel-Cantelli Lemma to conclude}
\begin{equation*}
    \limt \left[ E(U_t^2 \I(|U_t| \leq t^\delta ) | \fcalt) -  E(U^2_t|\fcalt) \right]  = 0\;\;\;\;a.s.
\end{equation*}
{\raggedright and apply the Cesàro's Mean Theorem to conclude}
\begin{scriptsize}
\begin{eqnarray*}
   \limn \frac{
                \sum\limits_{t=1}^n  E[U_t^2 \I(|U_t|\leq Z_t)|\fcalt ] -  E[U_t \I(|U_t|\leq Z_t)|\fcalt ]^2
              }{  
                  W_n^2
              } 
   &=& \limn \frac{
                \sum\limits_{t=1}^n  E[U_t^2|\fcalt] 
              }{  
                  n E\left( \D{l_t(\theta_r)}{i}^2  \right) 
              } \\
   &=&  1 \;\;\;\;a.s.
\end{eqnarray*}
\end{scriptsize}
{\raggedright By the generalized Chebyshev's inequality,}
\begin{eqnarray*}
 & & \hspace{-1cm} P\left[\left|E(U^2_t|\fcalt)-E(U^2_t I(|U_t| \leq t^\delta)|\fcalt)\right |>\varepsilon\right]  \\
 & & \hspace{1cm} \leq  \frac{1}{\varepsilon} E\left[ \left|E(U^2_t|\fcalt)-E(U^2_t I(|U_t| \leq t^\delta)|\fcalt) \right| \right] 
\end{eqnarray*}
{\raggedright and}
\begin{scriptsize}
\begin{eqnarray}
  \label{lem:derivada1a_assintoticaGARCHBEKK.e04}
  E\left[ \left|E(U^2_t|\fcalt)-E(U^2_t I(|U_t|\leq t^\delta)|\fcalt) \right| \right]  & & \notag \\
& & \hspace{-5cm}= E\left[ \left|E(U^2_t|\fcalt)-E(U^2_t I(|U_t|\leq t^\delta)|\fcalt) \pm E(U^2_t I(|U_t|>t^\delta)|\fcalt) \right| \right]  \notag\\
  & & \hspace{-5cm}=    E\left[ E(U^2_t I(|U_t|>t^\delta)|\fcalt)  \right]  \notag\\
  & & \hspace{-5cm}=    E\left[ U^2_t I(|U_t|>t^\delta) \right]  \notag\\
  & & \hspace{-5cm}\leq E\left[ U^4_t \right]^{1/2} E\left[I(|U_t|>t^\delta) \right]^{1/2}  \notag\\
  & & \hspace{-5cm}=    E\left[ U^4_t \right]^{1/2} \P(|U_t|>t^\delta)^{1/2} .
\end{eqnarray}
\end{scriptsize}
{\raggedright Using \eqref{lem:derivada1a_assintoticaGARCHBEKK.e03}, \eqref{lem:derivada1a_assintoticaGARCHBEKK.e05} and 
              \eqref{lem:derivada1a_assintoticaGARCHBEKK.e04}, for suitable $c>0$,}
\begin{eqnarray*}
  E\left[ \left|E(U^2_t|\fcalt)-E(U^2_t I(|U_t|\leq t^\delta)|\fcalt) \right| \right] 
  &\leq & c \frac{1}{ t^{\delta}} 
\end{eqnarray*}
{\raggedright which is $t$-summable.}

\item[(L3)] Using the notation $\Hdot:= \D{H_t}{i}$, 
 \begin{eqnarray}
   0 &\leq& E[U_t^4 \I(|U_t| \leq Z_t)| \fcalt] \notag \\
   &\leq&  E[U_t^4 | \fcalt] \notag \\
   &\leq&  E\left\{ 
                   \Tr \left[ X_t X'_t   H_t^{-1 } \Hdot H_t^{-1} 
                    -  H_t^{-1} \Hdot \right]^4
             | \fcalt  \right\} \notag \\
   &=&  E\left\{ 
                    \left[ \Tr \left( X_t X'_t   H_t^{-1 } \Hdot H_t^{-1}  \right)
                    -  \Tr\left(H_t^{-1} \Hdot\right) \right]^4
             | \fcalt  \right\} \notag \\
   &=&  E\left\{ 
                    \Tr \left( X_t X'_t   H_t^{-1 } \Hdot H_t^{-1}  \right)^4 \right.  \notag \\
  & &            -4 \Tr \left( X_t X'_t   H_t^{-1 } \Hdot H_t^{-1}  \right)^3 \Tr\left(H_t^{-1} \Hdot\right) \notag \\
  & &            +6 \Tr \left( X_t X'_t   H_t^{-1 } \Hdot H_t^{-1}  \right)^2 \Tr\left(H_t^{-1} \Hdot\right)^2 \notag \\
  & &            -4 \Tr \left( X_t X'_t   H_t^{-1 } \Hdot H_t^{-1}  \right) \Tr\left(H_t^{-1} \Hdot\right)^3 \notag \\
  & &  \left.      +\Tr\left(H_t^{-1} \Hdot\right)^4 | \fcalt  \right\} .
 \end{eqnarray}
{\raggedright Also, by Lemma \ref{lem:resultadosComteBEKKGARCH} (vi), for suitable $c\in(0,\infty)$,}
\begin{eqnarray}
       \left| \Tr \left( X_t X'_t   H_t^{-1 } \Hdot H_t^{-1}  \right)  \right |
   &=& \left| \Tr \left( H^{1/2}_t \varepsilon_t \left( H^{1/2}_t \varepsilon_t  \right)'    H_t^{-1 } \Hdot H_t^{-1}  \right) \right|  \notag \\
   &=& \left| \Tr \left( H^{1/2}_t \varepsilon_t \varepsilon'_t H_t^{1/2}    H_t^{-1 } \Hdot H_t^{-1}  \right)  \right| \notag \\
   &=& \left| \Tr \left( \varepsilon_t \varepsilon'_t H_t^{-1/2}  \Hdot H_t^{-1/2}  \right) \right|  \notag \\
   &\leq& \norm{\varepsilon_t \varepsilon'_t} \norm{ H_t^{-1/2}}^2  \norm{\Hdot }   \notag \\
   &\leq& c \norm{\varepsilon_t \varepsilon'_t}  \norm{\Hdot }  
\end{eqnarray}
{\raggedright and}
\begin{eqnarray}
       \left| \Tr \left( H_t^{-1/2} \Hdot  H_t^{-1/2} \right)  \right |
       &\leq&  c \norm{\Hdot} . 
\end{eqnarray}
{\raggedright We used the relation $|\Tr (ABC)| \leq \norm{A}\norm{B}\norm{C} $. Thus,}
\begin{scriptsize}
\begin{eqnarray}
  \label{lem:derivada1a_assintoticaGARCHBEKK.e05}
   E[U_t^4 \I(|U_t| \leq Z_t)| \fcalt] 
   &\leq&     c^4 E\left(\norm{\varepsilon_t\varepsilon'_t}^4 | \fcalt \right) \norm{\Hdot}^4  \notag \\
    & &      + 4 c^4 E\left(\norm{\varepsilon_t\varepsilon'_t}^3 | \fcalt \right) \norm{\Hdot}^4 \notag \\
    & &      + 6 c^4 E\left(\norm{\varepsilon_t\varepsilon'_t}^2 | \fcalt \right) \norm{\Hdot}^4 \notag \\
    & &      + 3 c^4 \norm{\Hdot}^4 \notag \\
    &\leq &    c_1 
                  E\left[ 
                    \norm{\varepsilon_t\varepsilon'_t}^2 
                    +\norm{\varepsilon_t\varepsilon'_t}^3 
                    +\norm{\varepsilon_t\varepsilon'_t}^4 
                    +1
                    | \fcalt 
                  \right]
                   \norm{\Hdot}^4. \notag \\
\end{eqnarray}
\end{scriptsize}
{\raggedright Adapting the proof of Lemma A.2 of Comte \& Lieberman \cite{comte_2003} and using B5, we conclude that}
\begin{equation}
 \label{eq:derivadaOrdem8Finita}
   E \sup\limits_{\theta\in\Theta_k} \left[ \norm{ \D{H_t}{i} (\theta) }^8 \right] < \infty.
\end{equation}
{\raggedright For $\delta_1\in(\frac{1}{2},1)$ and suitable $c_2,c_4\in(0,\infty)$, using the Chebyshev's and Jensen's inequalities, 
            \eqref{lem:derivada1a_assintoticaGARCHBEKK.e05} and \eqref{eq:derivadaOrdem8Finita},}
\begin{scriptsize}
\begin{eqnarray*}
  \P\left[  E[U_t^4 \I(|U_t| \leq Z_t)| \fcalt]  > t^{\delta_1} \right]   & & \\
     & & \hspace{-3cm} \leq  \P\left[  
                  E\left[ 
                    \norm{\varepsilon_t\varepsilon'_t}^2 
                    +\norm{\varepsilon_t\varepsilon'_t}^3 
                    +\norm{\varepsilon_t\varepsilon'_t}^4 
                    +1
                    | \fcalt 
                  \right]
                   \norm{\Hdot}^4 
        > \frac{ t^{\delta_1}}{c_1}  \right]   \\
     & & \hspace{-3cm}\leq  \frac {c^2_2} { t^{2 \delta_1 } } 
             E\left\{  
                  E\left[ 
                    \left(
                    \norm{\varepsilon_t\varepsilon'_t}^2 
                    +\norm{\varepsilon_t\varepsilon'_t}^3 
                    +\norm{\varepsilon_t\varepsilon'_t}^4 
                    +1
                    \right)
                   \norm{\Hdot}^4 
                    | \fcalt 
                  \right]^2
        \right\}   \\
     & & \hspace{-3cm} \leq  \frac {c^2_2} { t^{2 \delta_1 } } 
             E\left\{  
                  E\left[ 
                    \left(
                    \norm{\varepsilon_t\varepsilon'_t}^2 
                    +\norm{\varepsilon_t\varepsilon'_t}^3 
                    +\norm{\varepsilon_t\varepsilon'_t}^4 
                    +1
                    \right)^2
                   \norm{\Hdot}^8 
                    | \fcalt 
                  \right]
        \right\}   \\
     & & \hspace{-3cm}\leq  \frac {c^2_2} { t^{2 \delta_1 } } 
             E\left\{  
                    \left(
                    \norm{\varepsilon_t\varepsilon'_t}^2 
                    +\norm{\varepsilon_t\varepsilon'_t}^3 
                    +\norm{\varepsilon_t\varepsilon'_t}^4 
                    +1
                    \right)^2
             \right\}  
             E\left\{  
                   \norm{\Hdot}^8 
             \right\}   \\
    & & \hspace{-3cm}\leq \frac{c_3}{ t^{2 \delta_1 }} E\left[ \norm{\Hdot}^8 \right] 
\end{eqnarray*}
\end{scriptsize}
{\raggedright and }
\begin{eqnarray*}
  \sum\limits_{t=1}^\infty \P\left[   E[U_t^4 \I(|U_t| \leq Z_t)| \fcalt]  > t^{\delta_1} \right]  
  &\leq&  \sum\limits_{t=1}^\infty  \frac{c_3}{ t^{2 \delta_1 }} E\left[ \norm{\Hdot}^8 \right]  < \infty .
\end{eqnarray*}
{\raggedright By the Borel-Cantelli Lemma,}
\begin{equation*}
  \P\left[   E[U_t^4 \I(|U_t| \leq Z_t)| \fcalt]   > t^{\delta_1}\;\;\;\;i.o. \right]  =0 
\end{equation*}
{\raggedright and then}
 \begin{eqnarray*}
   \limn \sum\limits_{t=1}^n 
          \frac{
                E[U_t^4 \I(|U_t| \leq Z_t)| \fcalt]
              }{  
                  W_t^4 
              } 
    &\leq& \limn \sum\limits_{t=1}^n 
          \frac{
                1
              }{  
                t^{2-\delta_1}
              } 
    < \infty \;\;\;\;a.s.
 \end{eqnarray*}

\item[(L4)] $E\left(  \D{l_t(\theta_r)}{i}^2 \right)>0$, otherwise $\D{l_t(\theta_r)}{i} \equiv 0$ a.s. and $A_2$ would be $0$. Thus, 
                   using the stationarity,
\begin{scriptsize}
\begin{equation*}
   \limn \frac{
                W_n
              }{  
                  W_{n+1}
              } 
   = \limn \frac{
                 \left[ n E\left(  \D{l_t(\theta_r)}{i}^2  \right) \right]^{1/2}
              }{  
                 \left[ (n+1) E\left(  \D{l_t(\theta_r)}{i}^2  \right) \right]^{1/2}
              } 
              = 1 \;\;\;\;a.s.\;\;and\;\;\;\;
\end{equation*}
\begin{equation*}
    \limn W_n =  E\left(  \D{l_t(\theta_r)}{i}^2 \right)^{1/2} \limn \sqrt{n}  =  \infty.
\end{equation*}
\end{scriptsize}
\end{itemize}
\end{proof}
\end{lem}

\begin{lem}
 \label{lem:diversosParaH1BEKKGARCH}
 Let $\X=\left\{ X_t \right\}_{t\in \N}$ be a BEKK-GARCH(k),
 where conditions \eqref{eq:ergodicidadeBEKKGARCH} and B1-B5 are true. 
 Then,
\begin{equation*}
   E \sup\limits_{\theta\in\Theta_k} \left[ \left| \Tr(\Hdot(\theta) \Ht^{-1}(\theta) 
       - X_t X_t' \Ht^{-1}(\theta) \Hdot(\theta)\Ht^{-1}(\theta)) \right| \right] < \infty,
\end{equation*}
{\raggedright where $\Hdot:= D_\theta H_t$.}
\begin{proof}
{\raggedright Using Lemma \ref{lem:resultadosComteBEKKGARCH} and \eqref{eq:derivadaOrdem8Finita}, for a suitable $c\in(0,\infty)$,}
\begin{scriptsize}
\begin{eqnarray*}
  E\sup\limits_{\theta\in\Theta_k}\left[ \left| \Tr(\Hdot \Ht^{-1} - X_t X_t' \Ht^{-1} \Hdot\Ht^{-1}) \right| \right]
  &\leq& 
  E\sup\limits_{\theta\in\Theta_k} \left[ \norm{\Hdot} \norm{\Ht^{-1}} + \norm{X_t X_t'}\norm{ \Ht^{-1}}^2\norm{ \Hdot} \right] \\
  &\leq& 
  E\sup\limits_{\theta\in\Theta_k} \left[ 
  c
  \norm{\Hdot} + 
  c^2
  \norm{X_t X_t'}\norm{ \Hdot} \right] < \infty .
\end{eqnarray*}
\end{scriptsize}
\end{proof}
\end{lem}

\begin{lem}
 \label{thm:H1paraBEKKGARCH}
 Let $\X=\left\{ X_t \right\}_{t\in \N}$ be a BEKK-GARCH(r), of order $r$, 
 $k\not\geq r$,
 $\theta_r$ its true parameter,
$\thetahat_k$ the MLE of $\theta_r$.
 If conditions \eqref{eq:ergodicidadeBEKKGARCH} and B1-B5 are true, then
         \begin{equation*}
             \limn \frac{\loglnr{\thetahat_r} - \logl{\thetahat_k}}{n} > 0 \; a.s.
         \end{equation*}
\begin{proof}
Assuming
$p$ such as $p\geq k$ and $p\geq r$, 
\begin{scriptsize}
         \begin{equation*}
             \limn \frac{\loglnr{\thetahat_r} -  \logl{\thetahat_k}}{n}
             = 
             \limn \frac{\loglnr{\thetahat_r} - \loglnp{\thetahat_p} + \loglnp{\thetahat_p}  - \logl{\thetahat_k}}{n} .
         \end{equation*}
\end{scriptsize}
Applying Lemma \ref{lem:tetaImplicaLog} using the results above, we have
         \begin{eqnarray*}
             \limn \frac{\loglnr{\thetahat_r} - \loglnp{\thetahat_p} }{n}  = 0 \;\;\;\;a.s.
         \end{eqnarray*}
By \eqref{eq:logveroBEKKGARCH}, we see that we just need to prove
         \begin{eqnarray*}
             \limn \frac{\loglnp{\thetahat_p} - \logl{\thetahat_k} }{n}  > 0 \;\;\;\;a.s.
         \end{eqnarray*}
considering 
\begin{equation*}
  \logl{\theta} =  \sum\limitstUKbarN l_t(\theta)
\end{equation*}
and
\begin{equation*}
   l_t(\theta) 
        = \logpac{\frac{1}{\sqrt{(2\pi)^{m/2} \det(H_t)}} \exp\left( -\frac{1}{2 } x'_t H_t^{-1} x_t \right) }.
\end{equation*}
{\raggedright By Lemma \ref{lem:comportamento_assintotico_derivadasBEKKGARCH}, item (iii),}
\begin{equation*}
  E\left[ \left| l_t(\theta_r) \right| \right] < \infty
\end{equation*}
{\raggedright and then, using the SLLN,}
\begin{equation*}
 \limn \frac{\loglnp{\theta_r}}{ n} = \limn \frac{\sum\limitstUKbN  l_t(\theta_r) }{n } = E\left(   l_1(\theta_r) \right)  = c_1 < \infty \;\;\;\;a.s.
\end{equation*}
{\raggedright By the Mean Value Theorem, for $\dot{\theta} = s\theta_r + (1-s)\thetahat_p$, $s\in(0,1)$, 
              sufficiently large $n$ and $B_{\delta}(\theta_r)$ a sufficiently small neighborhood of $\theta_r$,}
\begin{eqnarray*}
 \left|  \frac{\sum\limitstUKbN  l_t(\thetahat_p)}{n } -  \frac{\sum\limitstUKbN  l_t(\theta_r)}{n } \right|
 &=& 
 \left|  \frac{\sum\limitstUKbN  D^1_{\theta} l_t(\dot{\theta})}{n } ( \thetahat_p - \theta_r) \right| \\
 &\leq& 
        \sup\limits_{\theta\in B_{\delta}(\theta_r)}  \norm{ \frac{\sum\limitstUKbN  D^1_{\theta} l_t(\dot{\theta})}{n } } \norm{( \thetahat_p - \theta_r)} \end{eqnarray*} 
{\raggedright Using the SLLN, Lemma \ref{lem:diversosParaH1BEKKGARCH} and the strong consistency of $\thetahat_p$,}
\begin{eqnarray}
\label{thm:H1paraBEKKGARCH.e01}
 \left|  \frac{\sum\limitstUKbN  l_t(\thetahat_p)}{n } -  \frac{\sum\limitstUKbN  l_t(\theta_r)}{n } \right|
 &\cas&  0 .
\end{eqnarray}
And thus,
\begin{eqnarray}
\limn \frac{\loglnp{\thetahat_p}}{n} = \limn \frac{\loglnp{\theta_r}}{n} = c_1.
\end{eqnarray}
Besides that, $\Theta_k\subset \Theta_{p}$ and $\thetahat_k$ is the MLE of $\theta_r$, thus
\begin{eqnarray*}
 \limn \frac{\logl{\thetahat_k}}{ n} 
 &\leq&\limn \frac{\loglnp{\thetahat_p}}{ n}  \\
 &=&  c_1 \;\;\;\;a.s.
\end{eqnarray*}
and then,
\begin{eqnarray*}
 \limn \frac{\logl{\thetahat_k}}{ n} 
 &\leq& \limsupn  \frac{\sum\limitstUKbN  l_t(\thetahat_k) }{n } = c_2 \leq c_1 \;\;\;\;a.s.
\end{eqnarray*}
Let $n_i$ be a subsequence of $n$ such as
\begin{eqnarray*}
 \limn  \frac{\sum\limits_{t=1+\kbar}^{n_i}  l_t(\thetahat_k) }{n_i } 
 &=& c_2 \;\;\;\;a.s.
\end{eqnarray*}
Using that $\Theta_k$ is compact, assume $n_j$ a subsequence of $n_i$ such as
\begin{eqnarray*}
\thetahat_k(n_j) \go \bar{\theta}_k\in\Theta_k \;\;\;\;a.s.
\end{eqnarray*}
And then,
\begin{eqnarray*}
  \limn \frac{\logl{\thetahat_k}}{n }
  &\leq& \limsupn \frac{\logl{\thetahat_k}}{n } \\
  &=& \lim\limits_{n_j \go \infty}  \frac{  \log L_{n_j} (  \thetahat_k(n_j) ) }{n_j} \;\;\;\;a.s.
\end{eqnarray*}
Applying the same argument used in  \eqref{thm:H1paraBEKKGARCH.e01}, we have
\begin{eqnarray*}
  \lim\limits_{n_i \go \infty}  \frac{  \log L_{n_j} (  \bar{\theta}_k ) }{n_j} 
  &=& E\left( l_1(\bar{\theta}) \right) \;\;\;\;and
\end{eqnarray*}
\begin{eqnarray*}
 \left|  \frac{\sum\limitstUKbN  l_t(\thetahat_k)}{n } -  \frac{\sum\limitstUKbN  l_t(\bar{\theta})}{n } \right|
 &\cas&  0 .
\end{eqnarray*}
Additionally,
\begin{eqnarray*}
 \limn \left[ \frac{\loglnp{\thetahat_p}}{ n} - \frac{\logl{\thetahat_k}}{ n} \right] 
  &\geq&  E\left[ \logpac{ \frac{ f(\bar{\theta}_k) }{ f(\theta_r)  }} \right] .\\
\end{eqnarray*}
By other hand,
$$
E\left[ \logpac{ \frac{ f(\bar{\theta}_k) }{ f(\theta_r)  }} \right]
$$
is the Kullback-Leibler divergence, which is positive if $ f(\bar{\theta}_k)\neq  f(\theta_r)$,
and, as $\theta_r\not\in \Theta_k\subseteq \R^{\gamma(k)}$, we have $\theta_r\neq \bar{\theta}_k$ and then, by B3,
$f(\bar{\theta}_k)\neq f(\theta_r)$.
Thus, we conclude
\begin{eqnarray*}
 \limn \left[ \frac{\logl{\thetahat_r}}{ n} - \frac{\logl{\thetahat_k}}{ n} \right] 
  &\geq&  E\left[ \logpac{ \frac{ f(\bar{\theta}_k) }{ f(\theta_r)  }} \right] = c > 0 .
\end{eqnarray*}
\end{proof}
\end{lem}

Lemmas \ref{lem:resultadosComteBEKKGARCH} and \ref{lem:comportamento_assintotico_derivadasBEKKGARCH} provides A1-A3 and Lemmas \ref{lem:derivada1a_assintoticaBEKKGARCH} and \ref{thm:H1paraBEKKGARCH} provides, respectively, A4 and A5. 
The EDC estimator's consistency is established using Theorem \ref{thm:consistencyEDC}.

\bibliography{referencias}

\end{document}